\documentclass[12pt]{article}
\usepackage{amsfonts}
\newtheorem{defi}{Definition}

\newtheorem{lm}{Lemma}
\newtheorem{thm}{Theorem}

\begin{document}

\title{On two-dimensional surface
attractors and repellers on
3-manifolds}
\author{V.~Grines\and
V.~Medvedev\and E.~Zhuzhoma}
\date{}
\maketitle \sloppy

\begin{abstract}

We show that if  $f: M^3\to M^3$  is an
$A$-diffeomorphism with a surface
two-dimensional attractor or repeller
$\mathcal{B}$ and $ M^2_ \mathcal{B}$ is a
supporting surface for $ \mathcal{B}$, then
$\mathcal{B} = M^2_ \mathcal{B}$ and there
is $k\geq 1$ such that:

1) $M^2_ \mathcal{B}$ is a union
$M^{2}_{1}\cup\dots\cup M^{2}_{k}$ of
disjoint tame surfaces such that every
$M^{2}_{i}$ is homeomorphic to the 2-torus
$T^2$.

2) the restriction of $f^k$ to $M^2_i$ $(i\in\{1,\dots,k\})$ is conjugate to Anosov
automorphism of $T^2$.

\end{abstract}

\section{Introduction}

One of the important question of Dynamical
Systems Theory is the relationship between a
fixed class of systems under consideration
and the topology of underlying manifolds.
This question is closely connected with a
structure of non-wandering set of a dynamic
system.   For example, Franks
\cite{Franks70} and Newhouse
\cite{Newhouse70} have  shown that any
codimension one Anosov diffeomorphism is
conjugate to a hyperbolic torus automorphism
(as a consequence, a manifold admitting such
diffeomorphisms is homeomorphic to the torus
$T^n$). A simple proof of this
Franks-Newhouse  theorem that uses foliation
theory techniques was obtained in
\cite{Hiraide2001}.

Recently  Grines and  Zhuzhoma
\cite{GrinesZh2005} proved that if a closed
$n$-manifold $M^n$, $n\geq 3$, admits a
structurally stable diffeomorphism  $f$ with
an orientable codimension one expanding
attractor, then $M^n$ is homotopy equivalent
to the $n$-torus $T^n$ and is homeomorphic
to $T^n$ for $n\neq 4$. Moreover, there are
no nontrivial basic sets different from the
codimension one expanding attractor in
nonwandering set of the   diffeomorphism
$f$. This allowed to them to classify, up to
conjugacy, structurally stable
diffeomorphisms having codimension one
expanding attractors and contracting
repellers on $T^n$.

A key point in the mentioned above results
is the existence of so-called hyperbolic
structure on a non-wandering set. More
precisely, let $ f: M\to M $ be a
diffeomorphism of a closed $m$-manifold $M$,
$m = \dim M\ge 2$, endowed with some
Riemannian metric $\rho$ (all definitions in
this section can be found in
\cite{KatokHassenblat95} and
\cite{Robinson-book99}, unless otherwise
indicated). Recall that a point $x\in M$ is
{\it non-wandering} if for any neighborhood
$U$ of $x$, $f^n(U)\cap U \neq \emptyset$
for infinitely many integers $n$. Then the
non-wandering set $NW(f)$, defined as the
set of all non-wandering points, is an
$f$-invariant and closed set. A closed
invariant set $\Lambda \subset M$ is {\it
hyperbolic} if there is a continuous
$f$-invariant splitting of the tangent
bundle $T_{\Lambda}M$ into stable and
unstable bundles $E^s_{\Lambda}\oplus
E^u_{\Lambda}$ with
$$ \Vert df^n(v)\Vert \le C\lambda ^n\Vert v\Vert ,\quad \Vert df^{-n}(w)\Vert
\le C\lambda ^n\Vert w\Vert , \quad \forall
v\in E^s_{\Lambda}, \forall w\in
 E^u_{\Lambda}, \forall n\in \mathbb{N},$$
for some fixed $C > 0$ and $\lambda < 1$.
For each $x\in \Lambda$, the sets $ W^s(x) =
\{y\in M: \lim\limits_{j\to \infty}\rho
(f^j(x),f^j(y))\to 0\} $, $ W^u(x) = \{y\in
M: \lim\limits_{j\to \infty} \rho
(f^{-j}(x),f^{-j}(y)) \to 0\} $ are smooth,
injective immersions of $E_x^s$ and $E_x^u$
that are tangent to $W_x^s$, $W_x^u$
respectively. $W^s(x)$, $W^u(x)$ are called
{\it stable} and {\it unstable manifolds} at
$x$.

An important class of dynamical systems is
made up of the diffeomorphisms satisfying
Smale's Axiom A \cite{Smale67}, the
so-called $A$-dif\-feo\-morp\-hisms. Given
such a diffeomorphism $f$, its recurrent
behavior is captured in its nonwandering set
$NW(f)$, which can be decomposed into
invariant topologically transitive pieces.
To be precise, a diffeomorphism $ f: M\to M
$ is an {\it $A$-dif\-feo\-morp\-hism} if
its non-wandering set $NW(f)$ is hyperbolic
and the periodic points are dense in
$NW(f)$. According to Smale's Spectral
Decomposition Theorem, $NW(f)$ is decomposed
into finitely many disjoint so-called basic
sets $\Omega_1$, $\ldots , \Omega_k$ such
that each $\Omega_i$ is closed,
$f$-invariant and contains a dense orbit
\cite{Smale67}. Following
\cite{AbrahamSmale70}, the pair $(a,b)$ is
said to be a {\it type} of basic set
$\Omega$ if
 $ a =  \dim E^s_x $
 and $ b =  \dim E^u_x $, where $ x\in \Omega $.

S.~Smale posed several kinds of basic sets:

(a) zero dimensional ones such as Smale's
horseshoe;

(b) one-dimensional ones such as so-called
Smale solenoids and its generalization
Smale-Williams solenoids (see
\cite{Williams67}, the name is suggested in
\cite{Pesin-book97});

(c) codimension one expanding attractors or
contracting repellers of DA-diffeomorphisms;

(d) basic sets of transitive Anosov
diffeomorphisms whose dimension equals to
the dimension of underlying manifold because
they coincide with the manifold.

It is well known that there is no
 restriction on the topology
 of underlying manifold when a basic set
 $\mathcal{B}$ is zero dimensional
 \cite{ShubSullivan75}.

But in the case of non zero dimensional the
situation is other in general (see above
mentioned examples). In addition,  we would
like to   notice the article \cite{Bothe81}
in which was proved that a closed orientable
$3$-manifold $M$ contains a knotted Smale
solenoid if and     only if
  $M$ has a lens space $L(p,q)$,
  with $p\ne 0,\pm1$, as a prime factor
  \footnote{Let $M$, $M_1$, $M_2$ are
  connected $3$-manifolds. We recall that $M$
   is a
connected sum of $M_1$ and $M_2$ and denote
this $M= M_1\sharp M_2$ if there are 3-cells
$B_i\subset M_i$ and embeddings $h_i : M_i -
Int B_i\to M$ ($i=1,2$) with $h_1(M_1 - Int
B_1)\cap h_2(M_2 - Int B_2) = h_1(\partial
B_1) = h_2(\partial B_2)$ and $M = h_1(M_1 -
Int B_1)\cup h_2(M_2 - Int B_2)$.  $M_1$,
$M_2$ are called factors. A connected
$3$-manifold $M$ is called prime manifold if
from condition $M= M_1\sharp M_2$ it follows
that exactly one manifold from $M_1$ and
$M_2$ must be prime. It is well known (see,
for example, \cite{Hempel-book}, theorem
3.15) that each compact 3-manifold can be
expressed as a connected sum of finite
number of prime factors. Let us recall that
a lens space is prime (see
\cite{Hempel-book}, ex. 3.12).}. This result
was repeated recently in
\cite{JiagNiWang2004}.

This paper is devoted to topological
classification of so-called surface basic
set of $A$-diffeomorphisms on smooth closed
orientable $3$-manifolds. Let us introduce a
concept of surface basic set.

\begin{defi} A basic set
$ \mathcal{B} $ of an A-diffeomorphism $f:
M^3\to M^3$ is called {\it surface basic
set} if $ \mathcal{B} $ belongs to an
$f$-invariant closed surface
$M^2_\mathcal{B} $ topologically embedded in
the 3-manifold $M^3$.
\end{defi}

The $f$-invariant surface $M^2_\mathcal{B} $
is called a {\it supporting surface} for $
\mathcal{B}$.

By definition, a supporting surface is not
necessary connected. But  it is obviously
that there is some power of diffeomorphism
$f$ for which every surface basic set has
connected supporting surface.

Let us recall that a basic set $\mathcal{B}$
of $A$-diffeomorphism $f:M\to M$ is called
an {\it attractor} if there is a closed
neighborhood $U$ of $\mathcal{B}$ such that
$f(U)\subset int~U$, $\bigcap\limits_{j\ge
0}f^j(U) = \mathcal{B}$.

If  $\mathcal{B}$ is a two-dimensional basic
set of $A$-diffeomorphism $f$ on a closed
3-manifold $M^3$ then, accordingly to
\cite{Plykin71} (theorem 3), $\mathcal{B}$
is either an attractor or repeller.

Recall that an  attractor is called an {\it
expanding attractor} if topological
dimension $\dim \mathcal{B}$ of
$\mathcal{B}$ is equal to the dimension
$\dim (E^u_{\mathcal{B}})$ of the unstable
bundle $ E^u_{\mathcal{B}} $ (the name is
suggested in \cite{Williams67},
\cite{Williams74}). A contracting repeller
of a diffeomorphism $f$ is an expanding
attractor for $f^{-1})$. Certainly, one can
consider a contracting repeller instead of
an expanding attractor, and vice versa. It
is well known that a codimension one
expanding attractor consists of the
$(\dim~M-1)$-dimensional unstable manifolds
$W^u(x)$, $x\in \mathcal{B}$, and is locally
homeomorphic to the product of $(\dim~M -
1)$-dimensional Euclidean space and a Cantor
set (see, for example, \cite{Plykin71},
theorem 2). The similar structure has
codimension one contracting repeller.
Therefore one use the notion {\it pseudotame
basic set}, meaning an expanding attractor
or contracting repeller.

As we know, the next Smale's question   is
open (see \cite{Smale67}, p. 785): is there
codimension one  basic set that is not
compact submanifolds and is not locally
product of Euclidean space and a Cantor set.

In section \ref{Proof of theorem 1} we shall
prove (in the  lemma \ref{type}) that
surface two-dimensional attractor (repeller)
$\mathcal{B}$ of an $A$-diffeomorphism
$f:M^3\to M^3$ has type $(2,1)$ ($(1,2)$).
It follows from there that  $\mathcal{B}$
is neither expanding attractor nor
contracting repeller.
 Moreover, we
will prove in section \ref{Proof of theorem
1}
 (lemma \ref{B=M}) that a two-dimensional
basic set $\mathcal{B}$ coincide with its
supporting surface $M^2_{\mathcal{B}}$.

Thus, Smale's question in the case under
consideration can be formulated as follows:
is there two-dimensional attractor
(repeller) which has  the type $(2,1)$
($(1,2)$) and  different from compact
submanifold?

Let us notice  that according to \cite{KMY}
there is an example of A-diffeomorphism of
closed three-manifold such that its
non-wandering set contains a two-dimensional
surface basic set whose supporting surface
is an essentially non-smoothly embedded
two-torus.

We notice also that according to
\cite{Pix77}, \cite{BoGr} there is  a
Morse-Smale diffeomorphism $f:S^3\to S^3$
with the $f$-invariant attracting surface
$S$ homeomorphic to the sphere $S^2$ that is
wildly embedded into $S^3$ (it is necessary
to emphasize that $S$ is not a basic set of
$f$ in this case).

However the first result of our paper claims
that a supporting surface for surface basic
set is a union of  tame tori.

Let us recall the concept of a  tame surface
embedded in $M^3$.

Let $D_0: \{(x,y,z)\in \mathbb{R}^3|
x^2+y^2\leq 1, z=0\}$ be  standard disk and
$M^2$ a surface embedded in  a
three-manifold $M^3$.

A surface $M^2$ is called {\it locally flat}
or {\it tame} if for any point $x\in M^2$
there is a neighborhood  $U_x$ of the point
$x$ in $M^3$ and homeomorphism $h_x:
\overline{U_x}\to \mathbb{R}^3$ such that
$h(\overline{U_x\cap M^2})=D_0$.

\begin{thm}
\label{M2-is-tame} Let $ f: M^3\to M^3 $ be
an $A$-diffeomorphism with the surface
two-dimensional basic set $ \mathcal{B} $
and  $ M^2_ \mathcal{B}$ is a supporting
surface for $ \mathcal{B}$. Then
$\mathcal{B} = M^2_ \mathcal{B}$ and there
is a number $k\geq 1$ such that $M^2_
\mathcal{B}$ is a union
$M^{2}_{1}\cup\dots\cup M^{2}_{k}$ of
disjoint tame surfaces such that every
$M^{2}_{i}$ is homeomorphic to the 2-torus
$T^2$.
\end{thm}

The next theorem explains  the dynamics of
restriction of the  diffeomorphism $f$ to a
surface basic set.

\begin{thm}
\label{f-anosov} Let the condition of
theorem \ref{M2-is-tame} are fulfilled, then
there is  number $k\geq 1$ such that the
restriction $f^k$ to $M^2_i$
$(i\in\{1,\dots,k\})$ is conjugate to Anosov
automorphism of $T^2$.
\end{thm}

{\it Acknowledgment}. The research was
partially supported by RFBR grant
02-01-00098.  This work was done while the
first and the third authors were visiting
Nantes University and Rennes 1 University
(IRMAR) respectively in March-June 2004.
They thank the support CNRS which made this
visits possible. They would like to thank
Francois Laudenbach, Andrey Pajitnov, Vadim
Kaimanovich, Anton Zorich,     for their
hospitality.

\section{Proof of theorem
\ref{M2-is-tame}} \label{Proof of theorem 1}

{\it Proof of theorem \ref{M2-is-tame}}
follows from the  next  lemmas \ref{type},
\ref{B=M} and \ref{tame} which will be
proved below.

\begin{lm}\label{type}
Let $\mathcal{B}$ be a two-dimensional
surface attractor (repeller) of
A-diffeomorphism $f$ of three-manifold
 $M^3$. Then  $\mathcal{B}$ has type
$(2,1)$ (resp. $(1,2)$).
\end{lm}
{\it Proof.} Suppose for definiteness that
$\mathcal{B}$ is a surface attractor and
$M^2_\mathcal{B}$ is the supporting surface
for $\mathcal{B}$. According to
\cite{Plykin71} (theorem 1), the unstable
manifold $W^{u}(x)$ belongs to $\mathcal{B}$
for any point $z\in\mathcal{B}$. Since the
restriction of $f$ to $\mathcal{B}$ is
transitive, $\dim W^u(z)$ does not depend on
the choice of $z\in \mathcal{B}$. Let us
notice that $\dim E^{u}_z = \dim W^{u}(z)$.
Thus it is sufficiently to  prove that
$\dim W^{u}(z)=1$.

Suppose the contrary. It follows from $\dim
W^u(z)\le\dim \mathcal{B}$ that either $\dim
W^u(z) = 0$ or $\dim W^u(z) = 2$.

If $\dim W^u(z) = 0$ then $\mathcal{B}$
would be an attracting orbit of the
diffeomorphism $f$ and, hence, a
zero-dimensional basic set. It contradicts
to supposition that $\dim \mathcal{B}=2$.

Suppose that  $\dim W^u(z) = 2$. As
$\mathcal{B}$ is nontrivial basic set, then
it contains infinite set
$Per(f)_{\mathcal{B}}$ of  periodic points
which are dense in $\mathcal{B}$. From
another hand, as $\mathcal{B}$ is a surface
basic set, then  for any point $p\in
Per(f)_{\mathcal{B}}$ the unstable manifold
$W^u(p)$ belongs to  the surface
$M^2_\mathcal{B}$ and consequently the set
$W^u(p)\setminus \{p\}$ does not contain
periodic points of the diffeomorphism $f$.
It contradicts to the fact that the set
$Per(f)_{\mathcal{B}}$ is dense in
$\mathcal{B}$. $\Box$

Let $\mathcal{B}$ be a two-dimensional
attractor of $A$-diffeomorphism $f: M^3\to
M^3$.

According to D.V. Anosov \cite{Anosov70} and
R. Bowen \cite{Bowen71} there is
 number
$k\geq 1$ such that the basic set
$\mathcal{B}$ can be represented as the
union of disjoint closed set $\mathcal{B}_1,
\dots, \mathcal{B}_k$ such that
$f(\mathcal{B}_i)=\mathcal{B}_{i+1}$
($\mathcal{B}_{k+1} =\mathcal{B}_{1}$) and
for any point $z\in \mathcal{B}_{i}$
$\overline{W^u(z)}=\mathcal{B}_{i}$. It
follows from  the proof of lemma \ref{type}
that $\dim W^u(z) = 1$ for any point $z\in
\mathcal{B}_i$ and $W^u(z)$ belongs to
$\mathcal{B}_i$.

Denote by $F^u_i$ the family  of unstable
manifolds  $W^u(z)$ for  all points $z\in
\mathcal{B}_i$.
\begin{lm}\label{B=M}
 Let $f: M^3\to M^3 $ be an
$A$-diffeomorphism whose non-wandering set
contains the surface two-dimensional
attractor
 $\mathcal{B}$ with the supporting
surface $M^2_\mathcal{B}$. Then:

1) $\mathcal{B} = M^2_\mathcal{B}$;

2)  the family $F^u_i$ is a continuous
foliation without singularities on
$M^{2}_{i}$.

3)  $M^2_\mathcal{B}$  is the union
$M^{2}_{1}\cup\dots\cup M^{2}_{k}$ of
disjoint surfaces such that each of them is
homeomorphic to the 2-torus $T^2$ and
$\mathcal{B}_i = M^2_i$.
\end{lm}

{\it Proof.}  Put $g=f^k$.  It is obviously
that for any $i\in \{1,\dots, k\}$ there is
a surface $M^2_i\subseteq M^2$ which is
supporting surface for the basic set
$\mathcal{B}_i$ of the diffeomorphism $g$.
Let us show that $\mathcal{B}_i=M^2_i$.

As by assumption $\dim \mathcal{B}_i = 2$
and $\mathcal{B}_i\subset M^2_i$, then
$\mathcal{B}_i$ contains a non-empty open
subset, say $V$, in the interior topology of
$M^2_i$  (see thm.~4.3
\cite{HurewiczWallmanm-book}).

Let $z$ be any point belonging to $int~V$.
According to \cite{Smale67}, there is
$\alpha>0$ such that the  point $z$  has a
closed neighborhood $U_z\subset V$ which is
homeomorphic to the direct product
$\check{W}^s_\alpha(z) \times
W^u_\alpha(z)$, where $\check{W}^s_\alpha(z)
= \mathcal{B}\cap W^s_\alpha(z)$ and
$W^s_\alpha(z)$, $W^u_\alpha(z)$ are a
closed  $\alpha$-neighborhoods of the point
$z$ in the some initial metric in the
manifold $W^s(z)$, $W^u(z)$ respectively. It
means that for any point $w\in U_z$ there is
a unique  pair of points   $w^s\in
\check{W}^s_\alpha(z)$, $w^u\in
W^u_{\alpha}(z)$ such that $w = W^u(w^s)\cap
W^s_\alpha(w^u)$. Let us define the
projection $\pi_z: U_z\to W^s(z)$ as
follows. For $w\in U_z$, put $\pi_z(w) =
w^s$.

Since $V$ is an open subset of $M^2$ and
$M^2$ is topologically embedded in $M^3$,
one can assume, without loss of generality,
that $U_z$ is homeomorphic to a closed
$2$-disk.

Put $U^s_z= \pi_z(U_z)$. Let us show that
there is closed  simple way $\gamma \subset
U^s_z$ such that $z \in int~\gamma$.

It is well known (see, for example,
\cite{Moise}) that any simple curve  is tame
on a disk. Thus there is an open
neighborhood $D_z\subset U_z$ of the point
$z$ such that:

1)  the set $W^u(z)\cap D_z$  consists of
exactly one connected component say $\lambda
\subset W^u(z)$;

2) the set $D\setminus \lambda$ consists of
two connected components say $D_1$, $D_2$.

3) $\pi_z(D_1)\cap \pi_z(D_2) = \emptyset$.

Let us choose any points $z_1 \in
\pi_z(D_1)$, $z_2 \in \pi_z(D_2)$. As the
sets $D_1\cup\lambda$, $D_2\cup\lambda$ are
path connected subsets and $\pi_z$ is
continuous map then $\pi_z(D_1\cup\lambda)$,
$\pi_z(D_2\cup\lambda)$ are also path
connected subsets. Moreover as
$\pi_z(D_1\cup\lambda)$,
$\pi_z(D_2\cup\lambda)$ are Hausdorff
subspaces (in induced topology), then there
is a simple arc $\gamma_1\subset
\pi(D_1\cup\lambda)$ joining the points
$z_1, z$ and  there is a simple arc
$\gamma_2\subset  \pi_z(D_2\cup\lambda)$
joining the points and $z_2, z$\footnote{See
for example \cite{Eng89}, proposition 6.3.12
(a) which claims that a Hausdorff space $X$
is pathwise connected if and only if for
every pair $x_1, x_2$ ($x_1\neq x_2$) there
exists a homeomorphic embedding  $h: I\to X$
of the closed unit interval in the space $X$
satisfying $h(0)=x_1, h(1)=x_2$ (e.i., $X$
is arcwise connected).}. Then the arc
$\gamma = \gamma_1\cup\gamma_2$ is a simple
arc joining the points $z_1$ and $z_2$. By
construction $z\in int~\gamma$.

Introduce a parameter $t\in [-1,1]$ on  the
arc $\gamma$ such that $\gamma (-1)=z_1$,
$\gamma (0)=z$ and $\gamma (1)=z_2$. For any
$t\in [-1,1]$ denote $l^u_t=W_\alpha^{
u}(\gamma(t))$ and put
$F^u_z=\bigcup\limits_{t\in [-1,1]}l^u_t$.

For any point $\gamma(t)$, the unstable
manifold $W_\alpha^u(\gamma(t))$ intersects
transversally the stable manifold
$W_\alpha^s(z)$ at  a unique point.

Let $V_0(z)$ be the disk   consisting of the
all points belonging to the union curves
$l^u_t$ ($t\in [-1,1]$). Due to the theorem
on continuous dependence of unstable
manifolds on initial conditions (see, for
example, \cite{Smale67}), the family of
curves $F^u_z$ is locally homeomorphic to a
family of straight lines parallel to the
axis Ox on Euclidean plain by means some
homeomorphism $h_z: V_0(z)\to \mathbb{R}^2$.
Hence, $F^u_z$ is a continuous foliation on
the closed disk $V_0(z)$.

Let $y\in \mathcal{B}_i$ be any point.   As
$\overline{W^{u}(y)}=\mathcal{B}_i$, then
there is a point $w\in int~\gamma$ such that
$[y,w]^u\subset W^{u}(y)$. Due to the
theorem on continuous dependence of unstable
manifolds on initial condition, there is an
open neighborhood $U_w\subset int~\gamma$ of
the point $w$  and an open neighborhood
$U_y\subset M^{2}_i$ of the point $y$ such
that for any point $y^{\prime}\in U_y$ there
is the point $w^{\prime} \in U_w$ such that
$[w^{\prime}, y^{\prime}]^u\subset
W^{u}(w^{\prime})$ and consequently the
point $y^\prime$ belongs to $\mathcal{B}_i$.
It means that the set $\mathcal{B}_i$ is
open. As $\mathcal{B}_i$ is closed then it
coincides with the surface $M^{2}_i$. As
$\mathcal{B}_i\cap \mathcal{B}_j =\emptyset$
for $i\neq j$ then also $M^{2}_i\cap M^{2}_j
=\emptyset$.

Thus  $\mathcal{B} = M^2_\mathcal{B}$ and
$M^2_\mathcal{B}$  is the union
$M^{2}_{1}\cup\dots\cup M^{2}_{k}$ of
disjoint surfaces.

It follows from above arguments that for any
point of $b\in \mathcal{B}_i$ there is a
neighborhood $U_b$  and homeomorphism $h_b:
U_b\to \mathbb{R}^2$ such that $h_b$ maps
the intersection of curves from $F^u_i$ with
$U_b$ on the family of straight lines which
are parallel to the axis $Ox$. Thus family
$F^u_i$ is continuous transitive foliation
without singularities on $M^{2}_i$.

Consequently  the surface $M^{2}_i$ is
homeomorphic to either the Klein bottle or
the torus. According to \cite{Kne} any
foliation without singularities on Klein
bottle must have at least one closed leaf.
As a consequence the  Klein bottle does not
admit transitive foliation. Thus the
manifold $M^{2}_i$ is homeomorphic to the
torus and the  lemma  is completely proved.
$\Box$

\begin{lm} \label{tame} The surface $M^2_i$  is tame.
\end{lm}{\it Proof.} According to the proof of
lemma \ref{B=M},  for any point $z\in M^2_i$
there exists a neighborhood $V_0(z)\subset
M^2_i$ such that:

1) $V_0(z)$ is homeomorphic to the direct
product $\gamma\times W^u_\alpha(z)$, where
$\gamma$ is simple curve belonging to
$W^s_\alpha(z)$;

2) $V_0(z)$ is the union of the smooth
curves $W^u_\alpha(w)$, $w\in \gamma$, which
belong to  leaves  of the foliation $F^u_i$
and form the foliation $\tilde{F}^u$ which
is given on the neighborhood $V_0(z)$;

3)  any curve $W^u_\alpha(w)$ intersects
$W^s_\alpha(z)$ in exactly one point;

4) the curve $\gamma$ is a local section (in
topological sense) for the foliation
$\tilde{F}^u$\footnote{We speak that
$\gamma$ is a local  section (in topological
sense) for the foliation $\tilde{F}^u$ if
for any point $d\in \gamma$ there is
neighborhood $V_d\subset V_0(z)$ such that
for any leaf $l\subset \tilde{F}^u$ with
$(l\cap V_d)\neq \emptyset$ the intersection
$(l\cap V_d)\cap \gamma$ consists of exactly
one point $x_l$ and $V_d\setminus\gamma$
consists  of two component each of which
contains exactly one component of the set
$l\setminus \{x_l\}$}.

Let us show that there are a neighborhood
$B_z$ of the point $z$ that is homeomorphic
to $3$-disk and  an embedding $h_z: B_z\to
\mathbb{R}^3$ such that $h_z(B_z\cap V_0(z))
= D_0$.

Without lost of generality we can suppose
that there exists  a neighborhood
$\tilde{B}_z$ (of the point $z$)  which is
 homeomorphic to
$3$-disk such that:

1) $W^s_{\alpha}(z)\subset \tilde{B}_z$,
$V_0(z)\subset \tilde{B}_z$;

2)  there is  diffeomorphism $g:\overline
{\tilde{B}_z}\to B_0^3$ such that
$g(W^s_\alpha(z))\subset Oxy$, $g(z) =
O(0,0,0)$, where $B_0^3$ is a
 closed unit ball in $\mathbb{R}^3$ .

Put $V= g(V_0(z))$, $Q= g(W^s_\alpha(z))$
and denote by $\hat{F}^u$ the foliation (on
$V$)  which is the image of the foliation
$\tilde{F}^u$ under the map $g$. By
construction the curve $\lambda = g(\gamma)$
is a local section (in topological sense)
for the foliation $\hat{F}^u$.  Let us
choose closed simple arc $\lambda_1\subset
int~\lambda$. Denote by $a$, $b$ the
endpoints of the arc $\lambda_1$ and denote
by $a_0$ ($b_0$) the point belonging to the
connected component of the set
$\lambda\setminus int~\lambda_1$ for which
the point $a$ ($b$) belongs to its boundary.
Denote by
 $\gamma_{a_0a}$ and $\gamma_{bb_0}$
 the closed simple ARCS belonging to
  $\lambda$ for which
 $a_0, a$ and $b,  b_0$ are the
 endpoints respectively.

Let us show that there exist the simple
piecewise linear arcs
$\tilde{\gamma}_{a_0a}$ and
$\tilde{\gamma}_{bb_0}$ belonging to $Q$
with endpoints $a_0, a$ and $b, b_0$
respectively  and satisfying the next
conditions:

1) $\tilde{\gamma}_{a_0a}\cap
\tilde{\gamma}_{bb_0}=\emptyset$;

2) $\tilde{\gamma}_{a_0a}\cap \lambda_1 =
a$, $\tilde{\gamma}_{bb_0}\cap \lambda_1=b$

3) the endpoints of each linear  link of the
curves $\tilde{\gamma}_{a_0a}$
($\tilde{\gamma}_{bb_0}$) belong to the
curve  $\gamma_{a_0a}$ ($\gamma_{bb_0}$).

Let us show the  existence of the curve
$\gamma_{a_0a}$  (the existence of the curve
$\gamma_{bb_0}$ may be shown similarly). Put
$\nu = \gamma_{a_0a}\setminus\{a\}$ and
choose an open neighborhood $U_\nu\subset V$
of the set $\nu$ such that

1) $U_\nu\cap(\lambda_1\cup
\gamma_{bb_0})=\emptyset$;

2)   $U_\nu$ admits a triangulation $\Sigma
= \bigcup\limits_{i\in\mathbb{Z}^+}\sigma_i$
such that for any point $x\in U_\nu$ any
neighborhood $U_x\subset U_\nu$ of the point
$x$  intersects only finite number of
simplexes from the union
$\bigcup\limits_{i\in\mathbb{Z}^+}
\sigma_i$.

Introduce a  parameter $t\in [0, \infty)$ on
the arc $\nu$ such that $\nu(0) = a_{0}$ and
$\nu(t)$ tends to $a$ as $t\to +\infty$.

Since $\nu(t)$ tends to $a$ as $t\in
+\infty$, there is the sequence of numbers
$0=t_{i_0}< t_{i_1}<,\dots,<  t_{i_k},
\dots$, where $t_{i_k}\to +\infty$ as $k\to
+\infty$, such that $\nu(t_{i_k}) \in
\sigma_{i_k}$,
$\sigma_{i_j}\cap\sigma_{i_{j+1}} \neq
\emptyset$, $int~\sigma_{i_j}\cap
int~\sigma_{i_{j+1}} = \emptyset$ and for
any $t>t_{i_k}$, $\nu(t)$ does not belong to
$\bigcup\limits_{j\leq k}\sigma_{i_j}$.
Denote by $l_k$ the piece of straight line
joining the points $\nu(t_{i_{k}})$ and
$\nu(t_{i_{k+1}})$, $k\in \mathbb{Z}^+$. By
construction, the sequence of the points
$\{\nu(t_{i_{k}})\}$ tends to $a$ as $k\to
+\infty$.  Then the set
$\tilde{\gamma}_{a_0a} =
\bigcup\limits_{i\in \mathbb{Z}^+}l_i\cup
\{a\}$ is a desired curve.

As  $\hat{F}^u$ is a continuous  foliation
consisting  of smooth curves which are
transversal to disk $Q$  there is a number
$N>0$ such that for any $c\in [-N,N]$ a
plain $P_c$ given by equation  $z = c$
intersects any leaf  of the foliation
$\hat{F}^u$ in exactly one point.

For any point $x\in \lambda$, denote by
$\tilde{L}^u_x$ the closed  arc such that:

1)  $\tilde{L}_x\subset L^u_x$, where
$L^u_x$ is the leaf of the foliation
$\hat{F}^u$ passing through the point $x$;

2) $\tilde{L}_x$ lyes between the plains
$P_{-N}: z=-N$, $P_{N}:z= N$.

Then the arc $\tilde{L}^u_{\nu(t_{i_{k}})}$
can be represented by equations:

$$x = x_k(z),  y = y_k(z), z=z, z\in
[-N,N].$$ Let us notice that by
construction, the point $\nu(t_{i_{k}})$ has
the coordinates $(x_k(0), y_k(0),0)$.

Denote by $S_k$ the disk  represented by the
next equations:

$x = x_k(z) + s_k(x_{k+1}(z)-x_k(z)), y =
y_k(z) + s_k(y_{k+1}(z)-y_k(z)), z=z,\\ z\in
[-N,N], s_k\in [0,1].$

Put $S_{a_0a}=\bigcup\limits_{k\in
\mathbb{Z}^+}S_k\cup \tilde{L}^u_a$.

 By construction, $S_{a_0a}$
is a piecewise smooth  disk. Using the
piecewise linear curve
$\tilde{\gamma}_{bb_0}$ we can construct a
piecewise smooth disk $S_{bb_0}$ which is
similar to $S_{a_0a}$.

Put $S_{ab}=
\bigcup\limits_{x\in\lambda_1}\tilde{L}^u_x$,
 $S=S_{a_0a}\cup S_{bb_0}\cup S_{ab}$,
$\nu_{-N} = S\cap P_{-N}$, $\nu_{N} = S\cap
P_{N}$.

It is not difficult to found the smoothly
embedded closed disks $S_1$,  $S_2$ such
that:

1) $S_1$, $S_2$  transversely intersect any
 plain $P_{c}$,  $c\in [-N,N]$;

2) $int~S_i\cap int~S =\emptyset$, $i=1,2$;

$S_i\cap S = \tilde{L}^u_{a_0}\cup
\tilde{L}^u_{b_0}$

3) $int~S_1\cap int~S_2 = \emptyset$,

$S_1\cap S_2 = \tilde{L}^u_{a_0}\cup
\tilde{L}^u_{b_0}$;

4) the boundary of the disk $S_i$ consists
of the curves $\tilde{L}^u_{a_0}$,
$\tilde{L}^u_{b_0}$ and the curves
$\nu^i_{-N} = S_i\cap P_{-N}$, $\nu^i_{N} =
S_i\cap P_{N}$, $i =1,2$.

5) the curves  $\nu^1_{-N}$, $\nu^2_{-N}$
($\nu^1_{N}$, $\nu^2_{N}$)
 form the boundary
of the closed disk $D_{-N}\subset P_{-N}$
($D_{N}\subset P_{N}$) which contains the
curve  $\nu_{-N}$ ($\nu_{N}$) dividing the
disk $D_{-N}$ ($D_{N}$) on two disks:
 $D^1_{-N}$ bounded by the
  curves $\nu_{-N}$, $\nu^1_{-N}$ and
  $D^2_{-N}$
bounded by the curves $\nu_{-N}$
   $\nu^2_{-N}$
 ($D^1_{N}$ bounded by the
  curves $\nu_{N}$, $\nu^1_{N}$ and
  $D^2_{N}$
bounded by the curves $\nu_{N}$
   $\nu^2_{N}$).

Denote  by $B_i$ ($i =1,2$) the closed set
bounded by the union $S\cup S_i\cup D^i_{-N}
\cup D^i_{N}$.

As intersection of the set  $B_i\cap P_c$ is
homeomorphic to the standard disk $D_0$,
then $B_i$ is homeomorphic to $D_0\times
[-N,N]$ and consequently $B_i$ is
homeomorphic
 to the standard ball $B_0:\{(x,y,z)\in \mathbb{R}^3|
x^2+y^2+z^2\leq 1\}$.

Then there is a neighborhood $U_S$ of the
disk $S$ such that $U_S\cap B_1$ ($U_S\cap
B_2$) is homeomorphic to $D_0\times [-1,0]$
($D_0\times [0,1]$). Thus there are a
neighborhood  $U_O$ of the point $O$ and
homeomorphism $h_O: \overline{U_O}\to
\mathbb{R}^3$ such that $h_O(U_O\cap
\hat{D})$ is the standard closed disk $D_0$.

The set   $B_z =g^{-1}(U_O)$ is a
neighborhood of the point $z$ in $M^3$ and
the map $h_z= h_O\circ g: B_z\to
\mathbb{R}^3$ satisfies to the following
condition: $h_z(B_z\cap V_0(z))$ is the
standard closed disk $D_0$. It means that
the surface $M^2_i$ is  tame. The lemma is
proved. $\Box$
\section{Proof of theorem \ref{f-anosov}}
\label{Proof of theorem 1}

{\it Proof of theorem \ref{f-anosov}}
follows from  lemmas \ref{hyper} and
\ref{hom} which will be proved below.

For any point $x\in M^2_i$ put
$L_i^s(x)=W^{s}(x)\cap M^2_i$ and denote
$F_i^s = \bigcup\limits_{x\in M^2_i}
L_i^s(x)$. It follows from the local
structure  of direct product and the proof
of lemma \ref{B=M} that there is $\alpha
>0$ such that for any point $z\in M^2_i$
there is the neighborhood $V_0(z)$ such
that:

1) $W_\alpha^u(z)\subset V_0(z)$ and for any
point $y\in W_\alpha^u(z)$ the  intersection
$W^s_\alpha(y)\cap V_0(z)$ consists  of a
simple curve $\gamma^s_y$;

2)  the family of curves
 $\bigcup\limits_{y\in W^u_\alpha(z)}\gamma^s_y$ is a
continuous   foliation on the neighborhood
$V_0(z)$, that is there is  a homeomorphism
$q_z: V_0(z)\to \mathbb{R}^2$ mapping the
family $\bigcup\limits_{y\in
W^u_{\alpha}(z)} \gamma^s_y$ to the set of
straight lines parallel to the axes $Ox$ on
the plain $\mathbb{R}^2$.

3) any curve $\gamma^s_y$ is locally section
  (in topological sense) to the
foliation $F^u_i$.

As $\mathcal{B}_i=M^2_i$ and for any point
$x\in \mathcal{B}_i$ the intersection
$W^s(x)\cap \mathcal{B}_i$ is dense in
$\mathcal{B}_i$, then any leaf  of the
foliation $F^s_i$ is dense in $M^2_i$. Thus
the foliation $F^s_i$ is  a transitive
foliation without singularities on the torus
$M^2_i$.

Let us represent the torus $M^2_i$ as the
factor space $\mathbb{R}^2/\Gamma$, where
$\Gamma$ is a discrete group of motions
$\gamma_{m,n}$ of the plain $\mathbb{R}^2$
given by the formulas $\gamma_{m,n}: \bar{x}
= x +m$, $\bar{y} = y +n$, $m, n\in
\mathbb{Z}$. Denote by $\pi:\mathbb{R}^2\to
M^2_i$ the natural projection  and $g_*$
automorphism of  the group $\Gamma$ induced
by diffeomorphism $g$ ($g =f^k$). Let us
notice that automorphism $g_*$ can be given
by the following way. Let
$\bar{g}:\mathbb{R}^2\to\mathbb{R}^2 $ be
the covering homeomorphism for $g:M^2_i\to
M^2_i$, that is $\pi\circ\bar{g}=g\circ\pi$,
and suppose that $\bar{g}$ is given by the
formulas: $\bar{g}: \bar{x} =g_1(x,y),
\bar{y} =g_2(x,y)$.   Then for any
$\gamma_{m,n}\in \Gamma$, one can put
$g_*(\gamma_{m,n}) =
\gamma_{m^\prime,n^\prime}$, where $m^\prime
= g_1(m,n)-g_1(0,0), n^\prime =
g_2(m,n)-g_2(0,0)$.

\begin{lm}\label{hyper}
The automorphism $g_*$ is hyperbolic, that
is the eigenvalues $\lambda_1$, $\lambda_2$
of the matrix
 $\mathbf{A}=
\left(\begin{array}{cc}
 a & b \\
 c & d
\end{array} \right)\in
SL(2,\mathbb{Z})$ which induces automorphism
$g_*$ satisfies to the condition
$|\lambda_1|<1$, $|\lambda_2|>1$.
\end{lm}

{\it Proof.} As  $F^{u}_i$ and  $F^{s}_i$
are transitive and transversal foliations on
the torus $M^2_i$  they form transitive
$2$-web on $M^2_i$. According to \cite{AGK}
(theorem 1), there are different irrational
numbers $\mu^u$, $\mu^s$ (which are Poincare
rotation numbers of the foliation $F^{u}_i$
and $F^{s}_i$ respectively) and
homeomorphism $\varphi: M^2_i\to M^2_i$ such
that $\varphi$ maps the foliations $F^{u}_i$
and $F^{s}_i$ to the linear foliation
$L_{\mu^u}$ and $L_{\mu^s}$ respectively
($L_{\mu^\sigma}$ is the image under the
projection $\pi$ of the   foliation
$\bar{L}_{\mu^\sigma}$ any leaf of whose is
given by the equation $y = \mu^\sigma x+c$,
$\sigma \in \{u,s\}, c\in \mathbb{R}^1$).

Denote by $\bar{F}^\sigma_i$ the foliation
on $\mathbb{R}^2$ which is covering for the
foliation $F^\sigma_i$ ($\sigma \in
\{s,u\}$). Then there is a homeomorphism
$\bar{\varphi}: \mathbb{R}^2\to
\mathbb{R}^2$ covering for the homeomorphism
$\varphi$ mapping the foliations
$\bar{F}^{u}_i$ and $\bar{F}^{s}_i$ to the
linear foliation $\bar{L}_{\mu^u}$ and
$\bar{L}_{\mu^s}$ respectively. It follows
from  there that  any leaf $l^u$ of the
foliation $\bar{F}^{u}_i$ intersects any
leaf $l^s$ of the foliation $\bar{F}^{s}_i$
at exactly one point.

Let us show now  that the automorphism $g_*$
is hyperbolic. Suppose the contrary. Let
$p\in Per(g)$ be a periodic points of a
period $m\geq 1$. Without lost of
generality, we can suppose  that $p =
\pi(O)$ (where $O$ is origin of the
coordinate system on  Euclidean plain).
Denote by $\bar{g}_m$ a covering
homeomorphism for $g^m$ such that
$\bar{g}_m(O)=O$. The set
$\mathcal{O}=\bigcup\limits_
{\gamma\in\Gamma}\gamma(O)$ is a lattice on
the plain $\mathbb{R}^2$. The matrix
$\mathbf{A}^m$ defines the automorphism
$g^m_*$. Denote by
$\mathcal{A}_m:\mathbb{R}^2\to \mathbb{R}^2$
the linear map determined by the matrix
$\mathbf{A}^m$. It follows from the
definition of $g_*$ that
$\bar{g}_m|_{\mathcal{O}}=
\mathcal{A}_m|_{\mathcal{O}}$. As a module
of the eigenvalues of the matrix $A^m$ is
equal to $1$ then  there  is a periodic
point $O_1\in \mathcal{O}$ ($O_1\neq O$) of
some period $l\geq 1$ of the map
$\mathcal{A}_m|_{\mathcal{O}}$. Consequently
$O_1$ is a periodic points of the
diffeomorphism  $\bar{g}_m$. Then the points
$O$, $O_1$ are fixed points of the
diffeomorphism $\bar{g}_m^{l}$.

Denote by
 $l^{u}_O$ ($l^{s}_{O_1}$)
the leaf of the foliation  $\bar{F}^u_i$
($\bar{F}^s_i$) passing through the point
 $O$ ($O_1$). As $l^{u}_O\cap l^{s}_{O_1} \neq \emptyset$ and
$l^{u}_O$, $l^{s}_{O_1}$ are invariant
unstable and stable manifold of the fixed
saddle (in topological sense) points $O$,
$O_1$ respectively then there are infinitely
many heteroclinic points of the
diffeomorphisms $\bar{g}_m^{l}$ belonging to
the intersection  $l^{u}_O \cap
l^{s}_{O_1}$. We get contradicts  with   the
fact that $l^{u}_O \cap l^{s}_{O_1} $
consists of  exactly one point. The lemma is
proved. $\Box$

  Denote by $\mathcal{G}$ the
linear automorphism of the torus $M^2_i$
such that $\mathcal{G}_*=g_*$.
 According to \cite{Franks70}
(proposition 2.1), there is a continuous
homotopic to identity map $h: M^2_i\to
M^2_i$ such that $\mathcal{G}h=hg$.

\begin{lm} \label{hom} The map $h$ is a
homeomorphism.
\end{lm}
{\it Proof.} Let
$\bar{h}:\mathbb{R}^2\to\mathbb{R}^2$ be a
covering map  for $h$. Let us divide the
proof of lemma into three steps.

{\it Step 1.} Let us show that if points
$\bar{x}, \bar{y}\in \mathbb{R}^2$
($\bar{x}\neq \bar{y}$) belong to the same
leaf $l^\sigma$ of  the foliation
$\bar{F}^\sigma_i$, then
$\bar{h}(\bar{x})\neq \bar{h}(\bar{y})$,
$\sigma\in \{s,u\}$.

Consider for definiteness the  case  $\sigma
= u$ (for $\sigma = s$, the proof is
similar) and suppose the contrary, that is
there are a leaf $l^u$ of the foliation
$\bar{F}^u_i$ and points $\bar{x},
\bar{y}\in l^u$ such that $\bar{h}(\bar{x})
= \bar{h}(\bar{y})$. Put $x=\pi(\bar{x})$,
$y=\pi(\bar{y})$, $L^u =\pi(l^u)$ and
$[x,y]^u\subset L^u$ is the closed arc with
endpoints  $x,y$. Let $p$ be any periodic
point of some period $l\geq 1$ of the
restriction of diffeomorphism $g$ to
$M^2_i$. Denote by $L^s_p$ a leaf of the
foliation $F^s_i$ passing through the point
$p$. We have two possibilities:

a)  the point $p$ belongs to $[x, y]^u$;

b)  the point $p$ does not belong to $[x,
y]^u$.

As the leaf $L^s_p$ is dense on the surface
$M_i^2$, then in   the case b) there is a
point $v\in L^s_p\cap (x, y)^u$.

 Consequently there are  two
cases:

$\bar{a})$ There is a point $\bar{p}\in
\pi^{-1}(p)$ belonging to the arc $[\bar{x},
\bar{y}]^u\subset l^u$;

$\bar{b})$ There are  a point $\bar{p}\in
\pi^{-1}(p)$ and a point $\bar{v}
\in\pi^{-1}(v)$ such that $\bar{p}$ and
$\bar{v}$ belong to the same leaf of the
foliation $\bar{F}^s$ and the point
$\bar{v}$ belong to the arc $(\bar{x},
\bar{y})^u\subset l^u$.

Let us consider a covering diffeomorphism
$\bar{g}_l$ for  the diffeomorphism $g^l$
such that $\bar{g}_l(\bar{p})=\bar{p}$. Then
in the both cases $\bar{a})$ and $\bar{b})$
we have $\rho(\bar{g}_l^n(\bar{x}),
\bar{g}_l^n(\bar{y}))\to +\infty$ as $n\to
+\infty$.

Let us notice that according to
\cite{Franks70} (lemma 3.4),  the map
$\bar{h}$ is proper\footnote{a map $\bar{h}$
is called proper if pre-image of a compact
set is a bounded set.}. Consequently,
according to \cite{Franks69},  there is a
number $r>0$ such that for any points
$\bar{x}_1, \bar{x}_2\in \mathbb{R}^2$
satisfying  a condition $\bar{h}(\bar{x}_1)=
\bar{h}(\bar{x}_2)$ an inequality
$\rho(\bar{x}_1, \bar{x}_2)< r$ is
fulfilled\footnote{For convenience we repeat
here arguments from \cite{Franks69}. Indeed,
as the map $\bar{h}$ is proper, then there
is $r>0$ such that pre-image of a
fundamental domain $\Pi$ of the group
$\Gamma$ belongs to the open disk
$D_{\frac{r}{2}} = \{(x,y)\in
\mathbb{R}|x^2+y^2< \frac{r^2}{4}\}$. Let
$\gamma \in \Gamma$ such that
$\gamma(\bar{h}(\bar{x}_1)) \in \Pi$. Then
as the map $h$ is homotopic to identity and
$\bar{h}(\bar{x}_1) = \bar{h}(\bar{x}_2)$ we
get $\gamma(\bar{h}(\bar{x}_1))=
\gamma(\bar{h}(\bar{x}_2))$,
$\bar{h}(\gamma(\bar{x}_1))=
\bar{h}(\gamma(\bar{x}_2))$. It follows from
there that $\rho(\gamma(\bar{x}_1),
\gamma(\bar{x}_2))< r$. As the map $\gamma$
is an isometry of $\mathbb{R}^2$ then
$\rho(\bar{x}_1, \bar{x}_2)< r$.}. As
$\mathcal{G}^l\circ h = h\circ g^l$ and $h$
is homotopic to identity then  there is a
covering  map $\bar{\mathcal{G}}_l$ for the
linear diffeomorphism $\mathcal{G}^l$ such
that $\bar{\mathcal{G}}_l\circ \bar{h} =
\bar{h}\circ g_l$. Then for any $n\in
\mathbb{Z}$ we have
$\bar{h}(\bar{g}^n(\bar{x}))=
\bar{\mathcal{G}}_l^n(\bar{h}(\bar{x}))=
\bar{\mathcal{G}}_l^n(\bar{h}(\bar{y}))=
\bar{h}(\bar{g}^n(\bar{y}))$. Consequently,
$\rho(\bar{g}^n(\bar{x}),
\bar{g}^n(\bar{y}))< r$. But it is
impossible as $\rho(\bar{g}_l^n(\bar{x}),
\bar{g}_l^n(\bar{y})\to +\infty$ as $n\to
+\infty$.

{\it  Step 2.} Let us show that for any
point $\bar{x}\in \mathbb{R}^2$ the
properties $\bar{h}(l^\sigma_{\bar{x}}) =
w^\sigma(\bar{h}(\bar{x}))$ is fulfilled,
where $l^\sigma_{\bar{x}}$ is the  leaf of
the foliation $\bar{F}_i^\sigma$ and
$w^\sigma(\bar{x})$ is the straight line
(passing through the point $\bar{x}$) which
is a pre-image   of invariant manifold
$W^{\sigma}(\pi(\bar{x}))$ of the linear
hyperbolic automorphism $\mathcal{G}$.

Consider for definiteness  a case   $\sigma
= s$ (in the case $\sigma = u$ the proof is
similar). First let us show that
$\bar{h}(l^\sigma_{\bar{x}})\subset
w^\sigma(\bar{h}(\bar{x}))$. Let $\bar{y}$
be any point belonging to
$l^\sigma_{\bar{x}}$ ($\bar{y}\neq
\bar{x}$). Put $x = \pi(\bar{x})$, $y =
\pi(\bar{y})$. As $\lim\limits_{n\to
+\infty}d(g^n(x),g^n(y))\to 0$ ($d$ is a
metric on the  torus $M^2_i$) then by
continuity of the map $h$ we have
$\lim\limits_{n\to
+\infty}d(h(g^n(x)),h(g^n(y)))=\lim\limits_{n\to
+\infty}d(\mathcal{G}^n(h(x)),
\mathcal{G}^n(h(y)))=0$. It follows from
there   that $h(y)\subset W^s(h(x))$. As
$\bar{h}$ is a covering map for $h$
 then $\bar{h}(l^s_{\bar{x}})\subset
 w^s(\bar{h}(\bar{x}))$.

 Let us show now that
 $\bar{h}(l^s_{\bar{x}}) =
  w^s(\bar{h}(\bar{x}))$.
Suppose the contrary. As
$\bar{h}(w^s(\bar{x}))$ is a connected set
which contains the point $\bar{h}(\bar{x})$
and  belongs to the straight line
$w^s(\bar{h}(\bar{x}))$,  then the image
under the map $\bar{h}$ of  (at least) one
component of the set $l^s_{\bar{x}}\setminus
\bar{x}$ is a bounded set on the straight
line $w^s(\bar{h}(\bar{x}))$.  It
contradicts to the fact that the map
$\bar{h}$ is proper.

{\it Step 3}. Let us show that the map $h$
is a homeomorphism. It follows from above
arguments:

1)  any point $\bar{x}$ of the plain
$\mathbb{R}^2$ is a unique point of the
intersections $l^s_{\bar{x}}\cap
l^u_{\bar{x}}$ and $w^s(\bar{x})\cap
w^u(\bar{x})$;

2)   the restriction of the map $\bar{h}$ to
each curve $l^s_{\bar{x}}$, $l^u_{\bar{x}}$
is a one-to-one  map on
$w^s(\bar{h}(\bar{x}))$,
$w^u(\bar{h}(\bar{x}))$ respectively.

It follows from 1) and 2)   that the map
$\bar{h}:\mathbb{R}^2\to \mathbb{R}^2$ is
one-to-one. Then the map $h:M^2_i\to M^2_i$
is a continuous one-to-one map and
consequently is a homeomorphism.

The lemma is completely proved. $\Box$

\bigskip

Dept. of Mathematics, Agriculture Academy of
Nizhny Novgorod, Nizhny  Novgorod, 603107
Russia

{\it E-mail address}: grines@vmk.unn.ru

{\it Current e-mail address}:
grines@math.univ-nantes.fr
\medskip

Dept. of Diff. Equat., Inst. of Appl.  Math.
and Cyber., Nizhny Novgorod,  State
University, Nizhny Novgorod, Russia

{\it E-mail address}: medvedev@uic.nnov.ru
\medskip

Dept. of Appl. Math., Nizhny Novgorod  State
Technical University, Nizhny  Novgorod,
Russia

{\it E-mail address}: zhuzhoma@mail.ru

{\it Current e-mail address}:
zhuzhoma@maths.univ-rennes1.fr

\end{document}